\documentclass[a4paper,11pt]{article}
\usepackage[T1]{fontenc}
\usepackage[cp1250]{inputenc}
\usepackage{lmodern}

\usepackage[english]{babel}
\usepackage{latexsym}
\usepackage{amsmath}
\usepackage{amsthm}
\usepackage{amssymb}
\usepackage{indentfirst}
\usepackage{amsfonts}
\usepackage{stackrel}
\usepackage{dsfont}
\usepackage{tikz}
\usepackage{slashbox}
\usepackage[mathscr]{eucal}
\usepackage{authblk}

\newtheorem{thm}{Theorem}
\newtheorem{lemma}[thm]{Lemma}
\newtheorem{cor}[thm]{Corollary}

\newenvironment{pro}{\begin{flushleft} \textbf{Proof}\\* \end{flushleft}}{\hfill\(\blacksquare\) \\ }

\newcommand{\reals}{\mathbb{R}}
\newcommand{\naturals}{\mathbb{N}}

\newcommand{\eps}{\varepsilon}

\newcommand{\Ffamily}{\mathcal{F}}

\newcommand{\Vfamily}{\mathcal{V}}
\newcommand{\Wfamily}{\mathcal{W}}

\newcommand{\C}{\mathscr{C}}

\newcommand{\Addresses}{{
  \bigskip
  \footnotesize
  Krukowski M. (corresponding author), \textsc{\L\'od\'z University of Technology, \ Institute of Mathematics, \ W\'ol\-cza\'n\-ska 215, \
90-924 \ \L\'od\'z, \ Poland}\par\nopagebreak
  \textit{E-mail address} : \texttt{krukowski.mateusz13@gmail.com}
}}

\begin{document}
\title{Short note on Smithson's paper concerning Arzel\`a-Ascoli theorem for multifunctions}
\author{Mateusz Krukowski}
\affil{\L\'od\'z University of Technology, Institute of Mathematics, \\ W\'ol\-cza\'n\-ska 215, \
90-924 \ \L\'od\'z, \ Poland}
\maketitle

\begin{abstract}
As indicated in the title, this paper was inspired by Smithson's work on Arzel\`a-Ascoli theorem. Up to this point, the author's interest in the topic has led him to generalize Arzel\`a-Ascoli theorem in the context of uniform spaces. This approach was primarily influenced by J.L. Kelley's \textit{General topology}. Soon after, the author learnt that in 1971, R.E. Smithson carried over Kelley's approach to the case of point-compact multifunctions (rather than single-valued functions). Similarly, in this paper the author would like to carry over the result, presented in his previous paper \textit{Arzel\`a-Ascoli theorem in uniform spaces}, to the case of point-compact multifunctions. 
\end{abstract}


\section{Introduction}

The author's interest in Arzel\`a-Ascoli theorem spiked nearly two years ago (in 2014) with the beginning of his PhD studies. In its basic formulation, the theorem provides the necessary and sufficient conditions under which every sequence of a given family of real-valued continuous functions, defined on $[0,1]$, has a uniformly convergent subsequence. These are equiboundedness and equicontinuity. 

Until now, the theorem has been generalized in multiple ways. For instance, Munkres (\cite{Munkres}, page 278) substitutes the interval $[0,1]$ with a compact space $X$. The space $C(X,\reals)$ is given the standard norm 
\begin{gather}\|f\| := \sup_{x\in X} |f(x)|.
\label{supnorm}
\end{gather}

\noindent
Still, equibounded and equicontinuous families coincide with the compact families in $C(X,\reals)$. If we allow for the change of topology determined by (\ref{supnorm}) to the topology of uniform convergence on compacta (which is not normable), we obtain a version of Arzel\`a-Ascoli theorem on $C(X,Y)$ for locally compact space $X$ and metric space $Y$ (\cite{Munkres}, page 290). In fact, metrizability of space $Y$ is not necessary and the version of Arzel\`a-Ascoli theorem for uniform spaces can be found in Kelley's (\cite{Kelley}, page 233) and Willard's monographs (\cite{Willard}, page 287). 

For most practical purposes (i.e. differential equations), the topology of uniform convergence on compacta is good enough to work with. Hence, nearly every extension of Arzel\`a-Ascoli theorem abandons the 'original' topology of uniform convergence in favour of the topology of uniform convergence on compacta. However, a nagging question remains: What if we do not want to resign from the uniform convergence? 

The first attempts to answer this question trace back to a paper by Bogdan Przeradzki and his study of differential equation of the form $x' = A(t)x + r(x,t)$ (comp. \cite{Przeradzki}). The paper gives a characterization of relatively compact subsets $\Ffamily \subset C^b(\reals,E)$, where $E$ is a Banach space. In addition to pointwise relative compactness and equicontinuity, the following condition was introduced:
\begin{description}
	\item[\hspace{0.4cm} (BP)] For any $\eps > 0$, there exist $T> 0$ and $\delta > 0$ such that if $\|f(T) - g(T)\| \leq \delta$ then $\|f(t) - g(t)\| \leq \eps$ for $t \geq T$ and if $\|f(-T) - g(-T)\|\leq \delta$ then $\|f(t) - g(t)\|\leq \eps$ for $t \leq -T$ where $f$ and $g$ are arbitrary functions in $\Ffamily$. 
\end{description}

\noindent
Robert Sta\'nczy (a student of Bogdan Przeradzki) further developed this approach, investigating Hammerstein equations in the space of bounded and continuous functions. In \cite{Stanczy}, he describes an application of this method to the Wiener-Hopf equations. However, as the approach became more abstract, the scope of possible applications quickly diminished. For instance in \cite{KrukowskiPrzeradzki}, the author together with Bogdan Przeradzki rewrote \textbf{(BP)} for $\sigma$-locally compact Hausdorff space $X$ and metric space $Y$. Tedious computations led to formulation of conditions under which the Hammerstein operator is compact.

Not daunted by the lack of natural applications, the author continued his study of Arzel\`a-Ascoli theorem. He introduced the quasimeasure of noncompactness and a new version of Darbo theorem (comp. \cite{KrukowskiDarbo}). Nevertheless, the author soon realized that (quasi)measures of noncompactness do not constitute a new technique in studying compactness. They merely provide a new language in which one can formulate old ideas. This is not meant as a derogatory remark. After all, shedding new light on the classical concepts is what makes them comprehensible.

The next chapter of study would not be possible without the influence of Grzegorz Andrzejczak and Bogdan Przeradzki. The first one introduced the author to the magical world of ultrafilters and Wallman compactification (\cite{GillmanJerison}, page 86). The second suggested using Gelfand-Naimark theorem (\cite{Kaniuth}, theorem 2.4.5 on page 69) as the convenient 'dictionary' between $C^b(T)$ and $C(\text{Wall}(T))$, where $T$ is a Tychonoff space. With this input, the author's work on the paper \textit{Arzel\`a-Ascoli theorem via Wallman compactification} (comp. \cite{KrukowskiWallman}) may be compared to that of a computer performing necessary calculations. The author owes a debt of gratitude to both his mentors.

The last chapter (as of 2016) of author's work was triggered by two monographs of Kelley and Willard (comp. \cite{Kelley} and \cite{Willard}). The author replaced the topology of uniform convergence on compacta with the topology of uniform convergence in the Arzel\`a-Ascoli theorem for $C(X,Y)$, where $X$ is a locally compact space and $Y$ is a Hausdorff uniform space. To this end, the author introduced (comp. \cite{Krukowskiuniform}) the compact extension property. Soon after submitting the paper to publication, the author realized that R.E. Smithson generalized the ideas of Kelley and Willard to the context of point-compact multifunctions (comp. \cite{Smithson}). The author's contribution in this paper is again the replacement the topology of uniform convergence on compacta with the topology of uniform convergence. This time, however, rather than Kelley's or Willard's monographs we use Smithson's paper as a reference point.

The main part of the paper is section \ref{sectionmainpart}, where we build on the Arzel\`a-Ascoli theorem for point-compact multifunctions (theorem \ref{arzela-ascoliforuniformconvergenceoncompacta} in \cite{Smithson}). The culminating point is the theorem \ref{ArzelaAscoli}, which characterizes compact subsets of continuous point-compact multifunctions, denoted by $\C(X,Y)$, with the topology of uniform convergence (rather than uniform convergence on compacta). The space $X$ is assumed to be locally compact and $Y$ is a Hausdorff uniform space.

In the spirit of \cite{FoxMorales}, \cite{LinRose} or \cite{Smithson}, we are not preoccupied with applications of our theorem. During his study, the author was struck by the beauty of Arzel\`a-Ascoli theorems and the accompanying concepts. Despite the lack of possible 'natural' applications, it is our firm belief that the presence of profound ideas (measures of noncompactness, Wallman compactification, Gelfand-Naimark theorem, uniform spaces to name a few) is a sufficient motivation for studying the topic.

\section{Preliminary notions}
\label{sectionrecap}

\subsection{Types of convergence}

Let $X$ be an arbitrary set and $(Y,\tau_Y)$ be a topological space. A \textit{multifunction} is a point to set correspondence $f : X \rightarrow Y$ such that for every $x \in X$, the set $f(x) \subset Y$ is nonempty (comp. \cite{FoxMorales}). A multifunction $f$ is said to be \textit{point-compact} if $f(x)$ is compact for every $x \in X$ (comp. \cite{Smithson}). In the sequel, we will only consider point-compact multifunctions. The set of all point-compact multifunctions is denoted by $Y^{mX}$. 

From this point onwards (unless otherwise stated), we assume that 
\begin{center}
\textit{$(Y,\Vfamily)$ is a Hausdorff uniform space.}
\end{center}

\noindent
Without loss of generality, we may assume that every set $V \in \Vfamily$ we are working with, is \textit{symmetric}. The \textit{uniformity of pointwise convergence} on $Y^{mX}$, which is denoted by $\Wfamily_{pc}$, is the uniformity, whose subbase sets are of the form
$$\bigg\{(f,g) \in Y^{mX} \times Y^{mX} \ : \ \forall_{\substack{y \in f(x)\\ z \in g(x)}}\ (y,g(x)) \cap V \neq \emptyset,\ (f(x),z) \cap V \neq \emptyset\bigg\},$$

\noindent
where $x \in X$ and $V \in \Vfamily$. Naturally, the \textit{topology of pointwise convergence}, denoted by $\tau_{pc}$, is the topology generated by the uniformity of pointwise convergence. However, in \cite{LinRose}, this topology is presented without resorting to the uniform structure. In order to do that, for $x \in X$ we define the $x$-projection $\pi_x : Y^{mX} \rightarrow X$ as the multifunction $\pi_x(f) = f(x)$. Now, the topology of pointwise convergence is defined as the topology having the open subbase of the forms 
\begin{gather}
\bigg\{f \in Y^{mX} \ : \ \pi_x(f) \cap U \neq \emptyset \bigg\} \hspace{0.2cm} \text{or} \hspace{0.2cm} \bigg\{f \in Y^{mX} \ : \ \pi_x(f) \subset U\bigg\}
\label{LinRosedefinition}
\end{gather}

\noindent
where $U \in \tau_Y$. 

In fact, both these descriptions are equivalent, as we shall shortly see (in theorem \ref{pointwiseconvergencecoincidence}). For convenience, we encompass two technical computations into the following lemma.

\begin{lemma}
For $f,g \in Y^{mX}$ and $x \in X$ we have 
$$\forall_{y \in f(x)}\ (y,g(x)) \cap V \neq \emptyset \ \Longleftrightarrow \ \pi_x(f) \subset V[g(x)]$$

\noindent
and for every $z \in g(x)$ we have 
$$(f(x),z) \cap V \neq \emptyset \ \Longleftrightarrow \ \pi_x(f) \cap V[z] \neq \emptyset.$$
\label{technicallemma}
\end{lemma}
\begin{pro}

For the first equivalence, we have 
\begin{gather*}
\pi_x(f) \subset V[g(x)] = \bigcup_{z \in g(x)}\ V[z] \ \Longleftrightarrow \ \forall_{y \in f(x)}\ \exists_{z \in g(x)}\ y \in V[z] \\
\Longleftrightarrow \ \forall_{y \in f(x)}\ \exists_{z \in g(x)}\ (y,z) \in V \ \Longleftrightarrow \ \forall_{y \in f(x)}\ (y,g(x)) \cap V \neq \emptyset.
\end{gather*}

\noindent
For the second equivalence, we have 
\begin{gather*}
\pi_x(f) \cap V[z] \neq \emptyset \ \Longleftrightarrow \ \exists_{y \in f(x)}\ y \in V[z] \\
\Longleftrightarrow \ \exists_{y \in f(x)}\ (y,z) \in V \ \Longleftrightarrow \ (f(x),z) \cap V \neq \emptyset 
\end{gather*}

\noindent
which ends the proof. 
\end{pro}

\begin{thm}
The topology of pointwise convergence and the one, whose subbase are the sets of the form \emph{(\ref{LinRosedefinition})}, coincide. 
\label{pointwiseconvergencecoincidence}
\end{thm}
\begin{pro}

The subbase set of the topology of pointwise convergence has the form 
\begin{gather*}
\bigg\{f \in Y^{mX} \ : \ \forall_{\substack{y \in f(x)\\ z \in g(x)}}\ (y,g(x)) \cap V \neq \emptyset\bigg\} \cap \bigg\{ f \in Y^{mX} \ : \ (f(x),z) \cap V \neq \emptyset\bigg\}\\
\stackrel{\text{lemma}\ \ref{technicallemma}}{=} \bigg\{f \in Y^{mX} \ : \ \pi_x(f) \subset V[g(x)]\bigg\} \cap \bigg\{ f \in Y^{mX} \ : \ \pi_x(f) \cap V[z] \neq \emptyset\bigg\}.
\end{gather*}

\noindent
It is now apparent that both topologies must agree, which ends the proof. 
\end{pro}

We recall the concept of \textit{uniformity of uniform convergence}, which is the main concern of Smithson's paper. For an extensive descriptions of this uniformity in the case of single-valued functions, we recommend Kelley's or Willard's monographs (\cite{Kelley} on page 226 or \cite{Willard} on page 280).

For a family $\Ffamily \subset Y^{mX}$ we define $\dagger : \Vfamily \rightarrow 2^{\Ffamily \times \Ffamily}$ by
\begin{gather}
\forall_{V \in \Vfamily} \ V^{\dagger} := \bigg\{ (f,g) \in \Ffamily \times \Ffamily \ : \ \forall_{\substack{x \in X\\ y \in f(x)\\ z \in g(x)}} \ (y,g(x)) \cap V \neq \emptyset,\ (f(x),z)\cap V \neq \emptyset\bigg\}.
\label{basetuc}
\end{gather}

\noindent
The family $\{V^{\dagger} \ : \ V \in \Vfamily\}$ is a base for uniformity $\Wfamily_{uc}$ on $\Ffamily$, which we call the \textit{uniformity of uniform convergence}. The precise verification of this statement would require checking all the conditions in theorem 2 in \cite{Kelley} on page 177. This is an easy task, so we skip the details.

Smithson elaborates on why this particular form of base sets (namely (\ref{basetuc})) is the most suitable. As the alternatives, he suggests the base sets of the form
$$\bigg\{(f,g) \in \Ffamily \times \Ffamily \ : \ f(x) \times g(x) \subset V\bigg\}$$

\noindent
and
$$\bigg\{(f,g) \in \Ffamily \times \Ffamily \ : \ f(x) \times g(x) \cap V \neq \emptyset \bigg\}.$$

\noindent
The first option leads to a uniformity, which is 'too little' while the second produces the uniformity which is 'too big'. Hence, despite being 'cumbersome' as Smithson put it, definition (\ref{basetuc}) is the most appropriate.

Again, as in the case of pointwise convergence, we have two descriptions of the topology of uniform convergence. Classicaly, it is induced by the uniformity of uniform convergence and has the base sets of the form  
\begin{gather}
V^{\dagger}[f] = \bigg\{g \in \Ffamily \ : \ \forall_{\substack{x \in X\\ y \in f(x)}} \ (y,g(x)) \cap V \neq \emptyset\bigg\} \cap \bigg\{g \in \Ffamily \ : \ \forall_{\substack{x \in X\\ z \in g(x)}}\ (f(x),z) \cap V \neq \emptyset\bigg\}
\label{tucopen}
\end{gather}

\noindent
where $V[y] = \{z \in Y \ : \ (y,z) \in V\}$ as in \cite{Kelley} on page 176. However, as in lemma \ref{technicallemma}, we may prove that $\tau_{uc}$ is generated by the base sets of the form 
$$\bigg\{f \in \Ffamily \ : \ \forall_{x \in X}\ \pi_x(f) \cap U \neq \emptyset \bigg\} \hspace{0.2cm} \text{or} \hspace{0.2cm} \bigg\{f \in \Ffamily \ : \ \forall_{x \in X}\ \pi_x(f) \subset U\bigg\}$$

\noindent
where $U \in \tau_Y$.

A concept, which is closely related to the uniformity of uniform convergence is the \textit{uniformity of uniform convergence on compacta}. Naturally, it is generated by the base sets of the form 
\begin{gather}
\bigg\{ (f,g) \in \Ffamily \times \Ffamily \ : \ \forall_{\substack{x \in D \\ y \in f(x)\\ z \in g(x)}} \ (y,g(x)) \cap V \neq \emptyset,\ (f(x),z)\cap V \neq \emptyset\bigg\}.
\label{basetucc}
\end{gather}

\noindent
where $V \in \Vfamily$ and $D \Subset X$ (meaning $D$ is a compact subset of $X$). The topology induced by the uniformity of uniform convergence on compacta is called the \textit{topology of uniform convergence on compacta}, which we denote by $\tau_{ucc}$. Repeating the reasoning that we have already done twice (for pointwise and uniform convergence), we conclude that base of this topology can be given by 
$$\bigg\{ g \in \Ffamily \ : \ \forall_{\substack{x \in D \\ y \in f(x)}} \ (y,g(x)) \cap V \neq \emptyset\bigg\} \cap \bigg\{ g \in \Ffamily \ : \ \forall_{\substack{x \in D \\ z \in g(x)}}\ (f(x),z)\cap V \neq \emptyset\bigg\}$$

\noindent
or, equivalently (again using lemma \ref{technicallemma}), by
$$\bigg\{f \in \Ffamily \ : \ \forall_{x \in D}\ \pi_x(f) \cap U \neq \emptyset \bigg\} \hspace{0.2cm} \text{or} \hspace{0.2cm} \bigg\{f \in \Ffamily \ : \ \forall_{x \in D}\ \pi_x(f) \subset U\bigg\}$$

\noindent
where $U \in \tau_Y$ and $D \Subset X$.

\subsection{Types of continuity}

Apart from the concepts of uniformities, we briefly review the notions of \textit{continuity} of multifunctions and, in particular, \textit{equicontinuity} (comp. \cite{FoxMorales}). A multifunction $f : X \rightarrow Y$ is said to be \textit{upper semi-continuous} if 
$$\forall_{W \in \tau_Y}\ f(x)\subset W \ \Longrightarrow \ \exists_{U_x \in \tau_X}\ f(U_x) \subset W.$$

\noindent
Furthermore, a multifunction is said to be \textit{lower semi-continuous} if 
$$\forall_{W \in \tau_Y}\ f(x) \cap W \neq \emptyset \ \Longrightarrow \ \exists_{U_x \in \tau_X}\ \forall_{y \in f(U_x)}\ f(y) \cap W \neq \emptyset.$$

\noindent
Not surprisingly, the function is called \textit{continuous} if it is both upper and lower semi-continuous. The set of continuous multifunctions, which are point-compact, will be denoted by $\C(X,Y)$. In the sequel, we will consider subfamilies of this space.

Another important concept is that of \textit{equicontinuity}. We say that the family $\Ffamily \subset Y^X$ is equicontinuous at $x \in X$ if for every $V \in \Vfamily$, there exists $U_x \in \tau_X$ such that
\begin{gather}
\forall_{f \in \Ffamily} \ f(U_x) \subset V[f(x)]
\label{equi1}
\end{gather}

\noindent
and
\begin{gather}
\forall_{\substack{z \in U_x\\ y \in f(x)}}\ f(z) \cap V[y] \neq \emptyset.
\label{equi2}
\end{gather}

\noindent
The family $\Ffamily$ is said to be equicontinuous if it is equicontinuous at every $x \in X$. The importance of equicontinuity is highlighted by the fact that if $\Ffamily$ is equicontinuous, then the topologies of pointwise convergence and the topology of uniform convergence on compacta coincide (comp. \cite{Smithson}).

\section{Arzel\`a-Ascoli theorems}
\label{sectionmainpart}

Our objective in this section is to present Arzel\`a-Ascoli theorem for the topology of uniform convergence in the setting of uniform spaces for point-compact multifunctions. From this this point onwards (unless otherwise stated), we assume that $X$ is not only a topological space, but also that
\begin{center}
\textit{$(X,\tau_X)$ is a locally compact space.}
\end{center}

\noindent
As the title suggests, the starting point for our considerations is theorem 9 in \cite{Smithson}.

\begin{thm}(Arzel\`a-Ascoli for uniform convergence on compacta)\\
Let $\C(X,Y)$ have the topology of uniform convergence on compacta $\tau_{ucc}$. A subfamily $\Ffamily \subset \C(X,Y)$ is relatively $\tau_{ucc}$-compact if and only if 
\begin{description}
	\item[\hspace{0.4cm} (AAucc1)] $\Ffamily$ is pointwise relatively compact, i.e. $\bigcup\ \{f(x) \ : \ f \in \Ffamily\}$ is relatively $\tau_Y$-compact for every $x \in X$, 
	\item[\hspace{0.4cm} (AAucc2)] $\Ffamily$ is equicontinuous.
\end{description}  
\label{arzela-ascoliforuniformconvergenceoncompacta}
\end{thm}

Our plan is to carry over the reasoning presented in \cite{Krukowskiuniform}, taking into account all the technicalities that arise with the change from single-valued functions to point-compact multifunctions. For $\Ffamily \subset Y^{mX}$ suppose that $\Wfamily_{pc}$ and $\Wfamily_{ucc}$ are uniformities of pointwise convergence and uniform convergence on compacta on $\Ffamily$, respectively. We recall (comp. \cite{Krukowskiuniform}) that $\Ffamily$ satisfies the \textit{finite extension property} if 
\begin{gather}
\forall_{V \in \Vfamily} \ \exists_{W \in \Wfamily_{pc}} \ W \subset V^{\dagger}
\label{finiteexconditionexplicitely}
\end{gather}

\noindent
and that it satisfies the \textit{compact extension property} if 
\begin{gather}
\forall_{V \in \Vfamily} \ \exists_{W \in \Wfamily_{ucc}} \ W \subset V^{\dagger}.
\label{compactexconditionexplicitely}
\end{gather}

Intuitively, condition (\ref{finiteexconditionexplicitely}) means that the topology $\tau_{pc}$ coincides with $\tau_{uc}$ on $\Ffamily$ (the inclusion $\tau_{pc} \subset \tau_{uc}$ always holds so (\ref{finiteexconditionexplicitely}) guarantees the reverse inclusion). Condition (\ref{compactexconditionexplicitely}) 'merely' implies that $\tau_{ucc}$ coincides with $\tau_{uc}$ on $\Ffamily$. Obviously if $\Ffamily$ satisfies finite extension property then it satisfies compact extension property.

\begin{thm}(Arzel\`a-Ascoli for uniform convergence)\\
Let $\C(X,Y)$ have the topology of uniform convergence $\tau_{uc}$. A subfamily $\Ffamily \subset \C(X,Y)$ is relatively $\tau_{uc}$-compact if and only if 
\begin{description}
	\item[\hspace{0.4cm} (AA1)] $\Ffamily$ is pointwise relatively compact and equicontinuous,
	\item[\hspace{0.4cm} (AA2)] $\Ffamily$ satisfies the finite extension property. 
\end{description}  
\label{ArzelaAscoli}
\end{thm}
\begin{pro}

An analogous, simple reasoning as in theorem 5 in \cite{Krukowskiuniform} shows that it suffices to prove that condition \textbf{(AA2)} is necessary. Suppose that \textbf{(AA2)} is not satisfied, which means that there exists $V \in \Vfamily$ such that $W\backslash V^{\dagger} \neq \emptyset$ for every $W \in \Wfamily_{pc}$. In what follows, $V^{\frac{1}{3}}$ will mean a symmetric set such that $V^{\frac{1}{3}} \circ V^{\frac{1}{3}} \circ V^{\frac{1}{3}} \subset V$ and likewise, we understand $V^{\frac{1}{9}}$ as a symmetric set such that $V^{\frac{1}{9}} \circ V^{\frac{1}{9}} \circ V^{\frac{1}{9}} \subset V^{\frac{1}{3}}$, where 
$$A \circ B = \bigg\{(a,b) \ : \ \exists_{c} \ (a,c) \in B, \ (c,b) \in A\bigg\}$$

\noindent
as in \cite{Kelley} on page 176. The existence of sets $V^{\frac{1}{3}}$ and $V^{\frac{1}{9}}$ follows from the axioms of uniformity (\cite{James} on page 103). In particular, our negation of \textbf{(AA2)} means that 
\begin{equation}
\begin{split}
&\bigg\{(f,g) \in \Ffamily\times\Ffamily \ : \ \forall_{\substack{x\in D\\ y \in f(x)\\ z \in g(x)}} \ (y,g(x)) \cap V^{\frac{1}{9}} \neq \emptyset,\ (f(x),z) \cap V^{\frac{1}{9}} \neq \emptyset \bigg\} \\
\cap \bigg\{(f,g) \in &\Ffamily\times\Ffamily \ : \ \exists_{\substack{x_{\ast} \in X\\ y_* \in f(x_*)}}\ (y_*,g(x_*)) \cap V = \emptyset \hspace{0.2cm} \text{or} \hspace{0.2cm} \exists_{\substack{x_* \in X\\ z_* \in g(x_*)}}\ (f(x_*),z_*) \cap V = \emptyset \bigg\} \neq \emptyset
\end{split}
\label{negofextcond}
\end{equation}

\noindent
for every finite set $D \subset X$. The family 
\begin{gather}
\bigg\{g \in C(X,Y) \ : \ \forall_{\substack{x \in X\\ y \in f(x)\\ z \in g(x)}}\ (y,g(x)) \cap V^{\frac{1}{9}}\neq \emptyset,\ (f(x),z)\cap V^{\frac{1}{9}} \neq \emptyset \bigg\}_{f \in \Ffamily}
\label{opencoverofF}
\end{gather}

\noindent
is a $\tau_{uc}$-open cover of $\overline{\Ffamily}$, the closure of $\Ffamily$ in $\tau_{uc}$. Indeed, if $\overline{f} \in \overline{\Ffamily}$, then by characterization of belonging to a closure, we have
$$(V^{\frac{1}{9}})^{\dagger}[\overline{f}] \cap \Ffamily \neq \emptyset \ \stackrel{(\ref{basetuc})}{\Longleftrightarrow} \ \exists_{f \in \Ffamily} \ \forall_{\substack{x\in X\\ y \in f(x)\\ z \in \overline{f}(x)}} \ (y,\overline{f}(x)) \cap V^{\frac{1}{9}} \neq \emptyset,\ (f(x),z) \cap V^{\frac{1}{9}} \neq \emptyset.$$

Since we assume that $\overline{\Ffamily}$ is $\tau_{uc}$-compact (and we aim to reach a contradiction), we can choose a finite subcover from (\ref{opencoverofF}), which means that there is a sequence $(f_k)_{k=1}^n \subset \Ffamily$ such that 
\begin{gather}
\forall_{g \in \overline{\Ffamily}} \ \exists_{k = 1,\ldots,n} \ \forall_{\substack{x \in X\\ y \in f_k(x)\\ z \in g(x)}} \ (y,g(x)) \cap V^{\frac{1}{9}} \neq \emptyset,\ (f_k(x),z) \cap V^{\frac{1}{9}} \neq \emptyset.
\label{whatsequence}
\end{gather}

Define, for every $k,l = 1,\ldots,n$, a set 
$$D_{k,l} := \bigg\{ x \in X \ : \ \exists_{y \in f_k(x)}\ (y,f_l(x)) \cap V^{\frac{1}{3}} = \emptyset \hspace{0.2cm} \text{or} \hspace{0.2cm} \exists_{z \in f_l(x)}\ (f_k(x),z) \cap V^{\frac{1}{3}} = \emptyset \bigg\}.$$

\noindent
Let $D$ be a set which consists of one element from each nonempty $D_{k,l}$. This set is finite and serves as a guard, watching whether each pair $f_k$ and $f_l$ 'drifts apart'. Its main task is the following implication
\begin{equation}
\begin{split}
&\forall_{\substack{x \in D\\ y \in f_k(x)\\ z \in f_l(x)}} \ (y,f_l(x)) \cap V^{\frac{1}{3}} \neq \emptyset,\ (f_k(x),z) \cap V^{\frac{1}{3}} \neq \emptyset \\ 
&\Longrightarrow \ \forall_{\substack{x \in X\\ y \in f_k(x)\\ z \in f_l(x)}} \ (y,f_l(x)) \cap V^{\frac{1}{3}} \neq \emptyset,\ (f_k(x),z) \cap V^{\frac{1}{3}} \neq \emptyset
\end{split}
\label{taskofD}
\end{equation}

\noindent
for every $k,l = 1,\ldots,n$. Indeed, suppose that for some $k,l = 1,\ldots,n$ there exist $x \in X,\ y\in f_k(x)$ such that $(y,f_l(x)) \cap V^{\frac{1}{3}} = \emptyset$. Then $D_{k,l} \neq \emptyset$, which implies $D \cap D_{k,l} \neq \emptyset$. We conclude that there exists $x_D \in D$ such that 
$$\exists_{y \in f_k(x_D)}\ (y,f_l(x)) \cap V^{\frac{1}{3}} = \emptyset \hspace{0.2cm} \text{or} \hspace{0.2cm} \exists_{z \in f_l(x_D)}\ (f_k(x_D),z) \cap V^{\frac{1}{3}} = \emptyset.$$ 

\noindent
An analogous reasoning works if there exist $x \in X,\ z \in f_l(x)$ such that $(f_k(x),z) \cap V^{\frac{1}{3}} = \emptyset$. This proves the implication (\ref{taskofD}).

By (\ref{negofextcond}), we pick $f_*,\ g_* \in \Ffamily$ such that 
\begin{gather}
\forall_{\substack{x\in D\\ y \in f_*(x)\\ z \in g_*(x)}} \ (y,g_*(x)) \cap V^{\frac{1}{9}} \neq \emptyset,\ (f_*(x),z) \cap V^{\frac{1}{9}} \neq \emptyset
\label{choiceoffgast1}
\end{gather}

\noindent
and
\begin{gather}
\exists_{\substack{x_* \in X\\ y_* \in f_*(x_*)}}\ (y_*,g_*(x_*)) \cap V = \emptyset \hspace{0.2cm} \text{or} \hspace{0.2cm} \exists_{\substack{x_* \in X\\ z_* \in g_*(x_*)}}\ (f_*(x_*),z_*) \cap V = \emptyset.
\label{choiceoffgast2}
\end{gather}

\noindent
Let $k_f, \ k_g \in \naturals$ be constants chosen as in (\ref{whatsequence}) for $f_{\ast}$ and $g_{\ast}$ respectively. For every $x \in D$, we have 
\begin{gather}
\left.\begin{array}{c}
\forall_{y \in f_{k_f}(x)}\ (y,f_*(x)) \cap V^{\frac{1}{9}} \neq \emptyset \\ 
\forall_{y \in f_*(x)}\ (y,g_*(x)) \cap V^{\frac{1}{9}} \neq \emptyset\\ 
\forall_{y \in g_*(x)}\ (y,f_{k_g}(x)) \cap V^{\frac{1}{9}} \neq \emptyset
\end{array}\right\} 
\Longrightarrow \ \forall_{y \in f_{k_f}} \ (y,f_{k_g}(x)) \cap V^{\frac{1}{3}} \neq \emptyset.
\label{bigimplication}
\end{gather}

\noindent
Indeed, take $x \in D$ and $y \in f_{k_f}(x)$. Then there exists $y' \in f_*(x)$ such that $(y,y') \in V^{\frac{1}{9}}$. Moreover, $(y',g_*(x))\cap V^{\frac{1}{9}} \neq \emptyset$ so there exists $y'' \in g_*(x)$ such that $(y',y'') \in V^{\frac{1}{9}}$. Finally, since $(y'',f_{k_g}(x)) \cap V^{\frac{1}{9}} \neq \emptyset$ we conclude that there exists $y''' \in f_{k_g}(x)$ such that $(y'',y''') \in V^{\frac{1}{9}}$. We have 
$$(y,y'),\ (y',y''),\ (y'',y''') \in V^{\frac{1}{9}} \ \Longrightarrow \ (y,y''') \in V^{\frac{1}{3}}$$

\noindent
so $(y,f_{k_g}(x)) \cap V^{\frac{1}{3}} \neq \emptyset$, which proves (\ref{bigimplication}). 

Similarly as above we prove that 
\begin{gather}
\forall_{\substack{x \in D\\ z \in f_{k_g}}} \ (f_{k_f}(x),z) \cap V^{\frac{1}{3}} \neq \emptyset.
\label{smallimplication}
\end{gather}

\noindent
By (\ref{taskofD}), we know that (\ref{bigimplication}) and (\ref{smallimplication}) imply
$$\forall_{\substack{x \in X\\ y \in f_{k_f}(x)\\ z \in f_{k_g}(x)}} \ (y,f_{k_g}(x)) \cap V^{\frac{1}{3}} \neq \emptyset,\ (f_{k_f}(x),z) \cap V^{\frac{1}{3}} \neq \emptyset.$$ 

\noindent
Finally (reasoning as above), for every $x \in X$ we have 
\begin{gather*}
\left.\begin{array}{c}
\forall_{y \in f_{\ast}(x)}\ (y,f_{k_f}(x)) \cap V^{\frac{1}{3}} \neq \emptyset \\
\forall_{y \in f_{k_f}(x)}\ (y,f_{k_g}) \cap V^{\frac{1}{3}} \neq \emptyset \\
\forall_{y \in f_{k_g}(x)}\ (y,g_*(x)) \cap V^{\frac{1}{3}} \neq \emptyset
\end{array}\right\} \ \Longrightarrow \ \forall_{y \in f_*(x)}\ (y,g_*(x)) \in V
\end{gather*}

\noindent
and
\begin{gather*}
\left.\begin{array}{c}
\forall_{z \in f_{k_f}(x)}\ (f_*(x),z) \cap V^{\frac{1}{3}} \neq \emptyset \\
\forall_{z \in f_{k_g}(x)}\ (f_{k_f}(x),z) \cap V^{\frac{1}{3}} \neq \emptyset \\
\forall_{z \in g_*(x)}\ (f_{k_g}(x),z) \cap V^{\frac{1}{3}} \neq \emptyset
\end{array}\right\} \ \Longrightarrow \ \forall_{z \in g_*(x)}\ (f_*(x),z) \in V
\end{gather*}

\noindent
which is a contradiction with (\ref{choiceoffgast2}). Hence, we conlcude that the finite extension property must hold. 
\end{pro}

The next corollary characterizes the relation between finite and compact extension property.

\begin{cor}
Let $\C(X,Y)$ have the topology of uniform convergence $\tau_{uc}$. For a subfamily $\Ffamily \subset \C(X,Y)$ satisfying \emph{\textbf{(AA1)}}, the following are equivalent: 
\begin{description}
	\item[\hspace{0.4cm} (C1)] $\Ffamily$ is relatively $\tau_{uc}$-compact. 
 	\item[\hspace{0.4cm} (C2)] $\Ffamily$ satisfies the finite extension property.
	\item[\hspace{0.4cm} (C3)] $\Ffamily$ satisfies the compact extension property. 
\end{description}  
\label{corollaryextprop}
\end{cor}
\begin{pro}

In theorem \ref{ArzelaAscoli} we proved the equivalence of \textbf{(C1)} and \textbf{(C2)}. Moreover, we already observed that finite extension property implies compact extension property (\textbf{(C2)} implies \textbf{(C3)}). Finally, if \textbf{(C3)} is satisfied, then the topologies $\tau_{ucc}$ and $\tau_{uc}$ conicide on $\Ffamily$ and realtive $\tau_{uc}$-compactness follows from theorem \ref{arzela-ascoliforuniformconvergenceoncompacta}.  
\end{pro}

\Addresses

\begin{thebibliography}{9}
\bibitem{Aliprantis}
	Aliprantis C.D., Border K.C. : \textit{Infinite Dimensional Analysis. A Hitchhiker's Guide}, Springer, 1999
\bibitem{Brezis}
	Brezis H. : \textit{Functional Analysis, Sobolev Spaces and Partial Differential Equations}, Springer, 2011
\bibitem{FoxMorales}
	Fox G., Morales P. : \textit{Non-Hausdorff multifunction generalization of the Kelley-Morse Ascoli theorem}, Pacific Journal of Mathematics, Vol. 64, No. 1 (1976)
\bibitem{GillmanJerison}
	Gillman L., Jerison M. : \textit{Rings of Continuous Functions}, D.Van Nostrand Company, 1960
\bibitem{James}
	James I.M. : \textit{Topologies and Uniformities}, Springer, 1999
\bibitem{Kaniuth}
	Kaniuth R. : \textit{A Course in Commutative Banach Algebras}, Springer, 2009
\bibitem{Kelley}
	Kelley J.L. : \textit{General Topology}, Springer, 1955
\bibitem{KrukowskiDarbo}
	Krukowski M. : \textit{Darbo-type theorem for quasimeasure of noncompactness}, arXiv: 1605.05229 (2016)
\bibitem{KrukowskiWallman}
	Krukowski M. : \textit{Arzel\`a-Ascoli theorem via Wallman compactification}, arXiv: 1602.05691 (2016)
\bibitem{Krukowskiuniform}
	Krukowski M. : \textit{Arzel\`a-Ascoli theorem in uniform spaces}, arXiv: 1602.05693 (2016)
\bibitem{KrukowskiPrzeradzki}
	Krukowski M., Przeradzki B. : \textit{Compactness result and its applications in integral equations}, arXiv:1505.02533 (2015)
\bibitem{LinRose}
	Lin Y.-F., Rose D.A. : \textit{Ascoli's theorem for spaces of multifunctions}, Pacific Journal of Mathematics, Vol. 34, No. 3 (1970)
\bibitem{Munkres}
  Munkres J. : \textit{Topology}, Prentice Hall, Inc., 2000
\bibitem{Przeradzki}
	Przeradzki B. : \textit{The existence of bounded solutions for differential equations in Hilbert spaces}, Annales Polonici Mathematici, LVI.2 (1992)
\bibitem{Smithson}
	Smithson R.E. : \textit{Uniform convergence for multifunctions}, Pacific Journal of Mathematics, Vol. 39, No. 1 (1971)
\bibitem{Stanczy}
	Sta\'nczy R. : \textit{Hammerstein equation with an integral over noncompact domain}, Annales Polonici Mathematici, 69 (1998)
\bibitem{Willard}
	Willard S. : \textit{General Topology}, Addison-Wesley Publishing Company, 1970
\end{thebibliography}
\end{document}